\documentclass[12pt]{amsart}

\newtheorem{theorem}{Theorem}[section]
\newtheorem{lemma}[theorem]{Lemma}
\newtheorem{prop}[theorem]{Proposition}
\newtheorem{cor}[theorem]{Corollary}

\theoremstyle{remark}
\newtheorem{remark}[theorem]{Remark}
\numberwithin{equation}{section}
\begin{document}
\title[Commutators of singular integrals]
{Commutators of singular integrals\\ on generalized $L^p$ spaces\\ with variable exponent}

\author{Alexei Yu. Karlovich}
\address{Departamento de Matem\'atica,
Instituto Superior T\'ecnico, Av. Rovisco Pais 1, 1049-001 Lisboa,
Portugal}
\email{akarlov@math.ist.utl.pt}

\author{Andrei K. Lerner}
\address{Department of Mathematics,
Bar-Ilan University, 52900 Ramat Gan, Israel}
\email{aklerner@netvision.net.il}

\begin{abstract}
A classical theorem of Coifman, Rochberg, and Weiss on
commutators of singular integrals is extended to
the case of generalized $L^p$ spaces with variable exponent.
\end{abstract}


\maketitle
\section{Introduction}
Let $T$ be a Calder\'on-Zygmund singular integral operator
$$
Tf(x):=\text{P.V.}\int_{\mathbb{R}^n}K(x-y)f(y)dy
$$
with kernel $K(x)=\Omega(x)/|x|^n$, where $\Omega$ is
homogeneous of degree zero, infinitely differentiable on the unit
sphere $S^{n-1}$, and $\int_{S^{n-1}}\Omega=0$.

All functions in the present paper are assumed to be real valued.
By $L^{\infty}_c$ we denote the class of all bounded functions on
$\mathbb{R}^n$ with compact support. Let $b$ be a locally
integrable function on $\mathbb{R}^n$. Consider the commutator
$[b,T]$ defined initially for any $f\in L^{\infty}_c$ by
$$[b,T]f:=bT(f)-T(bf).$$

Recall that the space $BMO(\mathbb{R}^n)$ consists of all
locally integrable functions $f$ such that
$$
\|f\|_{*}:=\sup_{Q}\frac{1}{|Q|}\int_Q|f(x)-f_Q|dx<\infty,
$$
where $f_Q:=|Q|^{-1}\int_Q f(y)dy$, the supremum is taken over all
cubes $Q\subset \mathbb{R}^n$ with sides parallel to the
coordinate axes, and $|Q|$ denotes the Lebesgue measure of $Q$.

A classical result of Coifman, Rochberg, and Weiss \cite{CRW76}
states that if $b\in BMO(\mathbb{R}^n)$, then $[b,T]$ is bounded on
$L^p(\mathbb{R}^n), 1<p<\infty$; conversely, if $[b,R_i]$ is
bounded on $L^p(\mathbb{R}^n)$ for every Riesz transform~$R_i$, then
$b\in BMO(\mathbb{R}^n)$. Janson \cite{Janson} observed that actually
for any singular integral $T$ (with kernel satisfying the
above-mentioned conditions) the boundedness of $[b,T]$ on $L^p(\mathbb{R}^n)$
implies $b\in BMO(\mathbb{R}^n)$.

An important role in the proving the latter implication is played by
a translation invariant argument, that is, by an obvious fact
that the translation operator is bounded on $L^p(\mathbb{R}^n)$.
Therefore, it is natural to ask whether an analogous result holds
if we replace $L^p(\mathbb{R}^n)$ by a more general function space for which the
continuity of translations may fail to hold.
We will consider the problem in \textit{generalized $L^p$ spaces
with variable exponent}.

Function spaces $L^{p(\cdot)}$ of Lebesgue type with variable exponent
$p$ were studied for the first time by Orlicz \cite{Orlicz31}.
Then Nakano considered spaces $L^{p(\cdot)}$ as an example
of his modular spaces \cite{Nakano50}. The theory of modular
spaces and, in particular, generalized Orlicz spaces generated by
Young functions with a parameter (Musielak-Orlicz spaces)
is documented in \cite{Musielak83}. The generalized $L^p$ spaces
with variable exponent are a special case of Musielak-Orlicz spaces.

Let $p:\mathbb{R}^n\to[1,\infty)$ be a measurable function. Consider the
convex modular (see \cite[Chap.~1]{Musielak83} for definitions and properties)
\[
m(f,p):=\int_{\mathbb{R}^n}|f(x)|^{p(x)}dx.
\]
Denote by $L^{p(\cdot)}(\mathbb{R}^n)$ the set of all Lebesgue
measurable functions $f$ on $\mathbb{R}^n$ such that
$m(\lambda f,p)<\infty$ for some $\lambda=\lambda(f)>0$.
This set becomes a Banach space with respect to the \textit{Luxemburg-Nakano norm}
\[
\|f\|_{L^{p(\cdot)}}:=\inf\Big\{\lambda>0: \ m(f/\lambda,p)\le 1\Big\}
\]
(see, e.g., \cite[Chap.~2]{Musielak83}).  Clearly, if $p(\cdot)=p$ is constant,
then the space $L^{p(\cdot)}(\mathbb{R}^n)$ is isometrically isomorphic
to the Lebesgue space $L^p(\mathbb{R}^n)$.

Observe, however, that spaces $L^{p(\cdot)}(\mathbb{R}^n)$ have attracted a great
attention only several years ago in connection with problems of
the boundedness of classical operators on $L^{p(\cdot)}(\mathbb{R}^n)$, which in
turn were motivated by some questions in fluid dynamics. We
mention here \cite{CUFN03,Diening02,DR,KS02,Nekvinda02,Ruzicka00} (see also
the references therein).

It is easy to see that $L^{p(\cdot)}(\mathbb{R}^n)$ fail to be
rearrangement-invariant, in general (see, e.g., \cite[Chap.~2]{BS}
for the definition and properties of rearrangement-invariant spaces).
This means that neither good-$\lambda$ technique nor rearrangement inequalities
may be applied for a generalization of some well-known results in
harmonic analysis to the case of $L^{p(\cdot)}(\mathbb{R}^n)$.
Also $L^{p(\cdot)}(\mathbb{R}^n)$ fail to be translation invariant,
in general (see \cite[Theorem~2.10]{KR91}).

If a measurable function $p:\mathbb{R}^n\to[1,\infty)$ satisfies
\begin{equation}\label{eq:reflexivity}
1<
p_-:=\operatornamewithlimits{ess\, inf}_{x\in\mathbb{R}^n} p(x),
\quad
\operatornamewithlimits{ess\, sup}_{x\in\mathbb{R}^n} p(x)=:p_+
<\infty,
\end{equation}
then the function
$$
p'(x):=p(x)/(p(x)-1).
$$
is well defined and satisfies (\ref{eq:reflexivity}) itself.

Denote by $\mathcal{M}(\mathbb{R}^n)$ the set of all
measurable functions $p:\mathbb{R}^n\to[1,\infty)$ such that
(\ref{eq:reflexivity}) holds and the Hardy-Littlewood maximal
operator $M$ is bounded on $L^{p(\cdot)}(\mathbb{R}^n)$.
Sufficient conditions guaranteeing $p\in\mathcal{M}(\mathbb{R}^n)$ are given in
\cite[Theorem~1.5]{CUFN03},
\cite[Theorem~3.5]{Diening02},
\cite[Theorem~2.14]{Nekvinda02} (see also \cite{KS02}
for weighted analogs).

Let $\mathcal{B}(X)$ be the class of all bounded sublinear
operators on a Banach lattice $X$ and let $\|A\|_{\mathcal{B}(X)}$
denote the operator norm of $A\in\mathcal{B}(X)$.

Our main result is the following.
\begin{theorem}\label{Main}
Suppose $p$ and $p'$ belong to $\mathcal{M}(\mathbb{R}^n)$.
\begin{enumerate}
\item[{\rm (a)}]
If $b\in BMO(\mathbb{R}^n)$, then $[b,T]$ is bounded on
$L^{p(\cdot)}(\mathbb{R}^n)$ and
$$
\|[b,T]\|_{\mathcal{B}(L^{p(\cdot)})}\le C_p\|b\|_*.
$$

\item[{\rm (b)}]
If $\Omega$ is odd, $b$ belongs to the Zygmund space $L\log L(Q)$
for every cube $Q\subset\mathbb{R}^n$
and $[b,T]$ is bounded on $L^{p(\cdot)}(\mathbb{R}^n)$,
then $b\in BMO(\mathbb{R}^n)$ and
$$
\|b\|_*\le C_p'\|[b,T]\|_{\mathcal{B}(L^{p(\cdot)})}.
$$
\end{enumerate}
\end{theorem}

Our proof of Part (a) is motivated by an analog of the
Fefferman-Stein theorem on the sharp maximal function for
$L^{p(\cdot)}(\mathbb{R}^n)$ proved recently by Diening and
R\r{u}\v zi\v cka \cite[Theorem~3.6]{DR}. To prove Part (a), we
combine a little bit more elaborate version of the latter result,
based on the so-called local sharp maximal function and on a
duality inequality due to the second author \cite[Theorem~1]{Le2},
with a sharp function inequality for commutators due to
Str\"omberg (see \cite{Janson}) and P\'erez
\cite[Lemma~3.1]{Perez}.

To prove Part (b), we first deduce that $[b,T]$ is also bounded
on $L^{p'(\cdot)}(\mathbb{R}^n)$. To make this step,
we have to pay by stronger requirements of oddness of the
kernel and of the local $L\log L$ integrability of $b$. Next,
using an interpolation argument, we conclude that $[b,T]$ is
bounded on $L^2(\mathbb{R}^n)$. This reduces the problem to the classical
situation.

We do not know whether assumptions on the kernel $K$ and on $b$ in
Part (b) of Theorem~\ref{Main} can be relaxed.

The paper is organized as follows. Section~\ref{sec:auxiliar} contains some
auxiliary results. In Section~\ref{sec:proof} we prove Theorem \ref{Main}.
Section~\ref{sec:concluding} contains some concluding remarks.
\section{Auxiliary results}
\label{sec:auxiliar}
\subsection{Duality and density in spaces $L^{p(\cdot)}(\mathbb{R}^n)$}
For $p$ satisfying (\ref{eq:reflexivity}) the function $p'$ is well defined
and one can equip the space
$L^{p(\cdot)}(\mathbb{R}^n)$ with the \textit{Orlicz type norm}
$$
\|f\|_{L^{p(\cdot)}}^0
:=
\sup\left\{\int_{\mathbb{R}^n}|f(x)g(x)|dx\ :\
g\in L^{p'(\cdot)}(\mathbb{R}^n),
\quad
\|g\|_{L^{p'(\cdot)}}\le 1 \right\}.
$$
This norm is equivalent to the Luxemburg-Nakano norm (see \cite[Theorem~2.3]{KR91}):
\begin{equation}\label{eq:Orlicz-Luxemburg}
\|f\|_{L^{p(\cdot)}}
\le
\|f\|_{L^{p(\cdot)}}^0
\le
r_p\|f\|_{L^{p(\cdot)}}
\quad \Big(f\in L^{p(\cdot)}(\mathbb{R}^n)\Big),
\end{equation}
where
$$
r_p:=1+1/p_--1/p_+.
$$
\begin{lemma}\label{le:Hoelder}
{\rm (see \cite[Theorem~2.1]{KR91}).}
Let $p:\mathbb{R}^n\to[1,\infty)$ be a measurable function
satisfying {\rm (\ref{eq:reflexivity})}. If $f\in L^{p(\cdot)}(\mathbb{R}^n)$
and $g\in L^{p'(\cdot)}(\mathbb{R}^n)$, then $fg$ is integrable on
$\mathbb{R}^n$ and
$$
\int_{\mathbb{R}^n}|f(x)g(x)|dx
\le
r_p\|f\|_{L^{p(\cdot)}}
\|g\|_{L^{p'(\cdot)}}.
$$
\end{lemma}

   From \cite[Theorem~2.11]{KR91} we get the following.
\begin{lemma}\label{le:density}
Let $p:\mathbb{R}^n\to[1,\infty)$ be a measurable function
satisfying {\rm (\ref{eq:reflexivity})}. Then $L_c^\infty$
is dense in $L^{p(\cdot)}(\mathbb{R}^n)$ and in
$L^{p'(\cdot)}(\mathbb{R}^n)$.
\end{lemma}
\subsection{Pointwise estimates for sharp maximal functions}
Given $f\in L_{\rm loc}^1(\mathbb{R}^n)$, the
Hardy-Littlewood maximal function is defined by
$$
Mf(x):=\sup_{Q\ni x}\frac{1}{|Q|}\int_{Q}|f(y)|dy.
$$
For $\delta>0$ and $f\in L^\delta_{\rm loc}(\mathbb{R}^n)$, set also
$$
f_{\delta}^{\#}(x):=\sup_{Q\ni x}\inf_{c\in\mathbb{R}}
\left(\frac{1}{|Q|}\int_Q|f(y)-c|^{\delta}dy\right)^{1/\delta}.
$$

The non-increasing rearrangement (see, e.g., \cite[Chap.~2,
Section~1]{BS}) of a measurable function $f$ on $\mathbb{R}^n$ is
defined by
$$
f^*(t):=
\inf\Big\{\lambda>0:|\{x\in\mathbb{R}^n:|f(x)|>\lambda\}|\le t\Big\}
\quad(0<t<\infty).
$$
Set also $f^{**}(t):=t^{-1}\int_0^tf^*(\tau)d\tau$.

For a fixed $\lambda\in(0,1)$ and a given measurable function $f$
on $\mathbb{R}^n$, consider the local sharp maximal function
$M^{\#}_{\lambda}f$ defined by
$$
M^{\#}_{\lambda}f(x) := \sup_{Q\ni x}\inf_{c\in\mathbb{R}}
\big((f-c)\chi_Q\big)^*\left(\lambda|Q|\right).
$$
In all above definitions the supremums are taken over all cubes
$Q\subset\mathbb{R}^n$ containing $x$.
\begin{prop}\label{pr:relation}
If $\delta>0,\lambda\in(0,1)$, and $f\in L_{\rm loc}^\delta(\mathbb{R}^n)$,
then
\begin{equation}\label{eq:relation-1}
M_\lambda^\# f(x)\le (1/\lambda)^{1/\delta}f_\delta^\#(x)
\quad (x\in\mathbb{R}^n).
\end{equation}
\end{prop}
\begin{proof}
Let $\varphi\in L_{\rm loc}^\delta(\mathbb{R}^n)$ and $x\in\mathbb{R}^n$.
For every cube $Q$ containing $x\in\mathbb{R}^n$, by Chebyshev's
inequality,
\begin{equation}\label{eq:relation-2}
\big(|\varphi|^\delta\chi_Q\big)^*(\lambda|Q|) \le
\frac{1}{\lambda|Q|}\int_Q|\varphi(y)|^\delta dy.
\end{equation}
On the other hand, in view of \cite[Chap.~2, Proposition~1.7]{BS},
\begin{equation}\label{eq:relation-3}
\big(|\varphi|^\delta\chi_Q\big)^*=\big[(\varphi\chi_Q)^*\big]^\delta.
\end{equation}
Take $\varphi=f-c$ with $c\in\mathbb{R}$. Then from (\ref{eq:relation-2})
and (\ref{eq:relation-3}) we get
$$
\big((f-c)\chi_Q\big)^*(\lambda|Q|)
\le
(1/\lambda)^{1/\delta}
\left(\frac{1}{|Q|}\int_Q|f(y)-c|^\delta dy\right)^{1/\delta}.
$$
Taking the infimum over $c\in\mathbb{R}$ and then the supremum
over all cubes $Q\subset\mathbb{R}^n$ containing $x$, we obtain
(\ref{eq:relation-1}).
\end{proof}
\begin{theorem}\label{th:AlvarezPerez}
{\rm (see \cite[Theorem~2.1]{AP}).}
If $0<\delta<1$, then for every $f\in L^{\infty}_c$,
$$
(Tf)_\delta^\#(x)\le c_{\delta,n}Mf(x)
\quad (x\in\mathbb{R}^n).
$$
\end{theorem}

A sharp function inequality for the commutator $[b,T]f$ was proved
by Str\"omberg (the proof is contained in \cite{Janson}).
Afterwards, a more precise version of this result was given by
P\'erez \cite{Perez}. We will need the following corollary from
\cite[Lemma~3.1]{Perez}.
\begin{theorem}\label{th:Perez}
If $0<\delta<1$, then for every $b\in BMO$ and $f\in L_c^\infty$,
$$
([b,T]f)_{\delta}^{\#}(x)
\le
c_{\delta,n}\|b\|_*\big(M(Tf)(x)+MMf(x)\big)
\quad (x\in\mathbb{R}^n).
$$
\end{theorem}

Actually, Theorems~\ref{th:AlvarezPerez} and \ref{th:Perez} were
proved in \cite{AP} and \cite{Perez}, respectively, for smooth
functions, but exactly the same proofs work for $L^{\infty}_c$
functions as well.
\begin{theorem}\label{th:Lerner}
{\rm (see \cite[Theorem~1]{Le2}).} For a function $g\in L_{\rm
loc}^1(\mathbb{R}^n)$ and a measurable function $\varphi$
satisfying
\begin{equation}\label{eq:Lerner}
|\{x:|\varphi(x)|>\alpha\}|<\infty \quad\mbox{for all}\quad
\alpha>0,
\end{equation}
one has
$$
\int_{\mathbb{R}^n}|\varphi(x)g(x)|dx \le
c_n\int_{\mathbb{R}^n}M_{\lambda_n}^\# \varphi(x)Mg(x)dx.
$$
\end{theorem}
\subsection{On the boundedness of singular integral operators}
\begin{theorem}\label{th:CZ}
If $p,p'\in\mathcal{M}(\mathbb{R}^n)$, then there exists a
constant $c_p$ such that for any $f\in
L^{p(\cdot)}(\mathbb{R}^n)$,
$$
\|Tf\|_{L^{p(\cdot)}}\le c_p\|f\|_{L^{p(\cdot)}}.
$$
\end{theorem}
\begin{remark}
This result was proved by Diening and R\r{u}\v zi\v cka
\cite[Theorem~4.8]{DR} under the assumptions
$p\in\mathcal{M}(\mathbb{R}^n)$ and
$(p/s)'\in\mathcal{M}(\mathbb{R}^n)$ for some $s\in(0,1)$. Their
proof is based on an analog of the Fefferman-Stein theorem on the
sharp maximal function proved in the same paper
\cite[Theorem~3.6]{DR}. We shall give a little bit different proof
for the sake of completeness. See also Remark~\ref{rem:final}
below.
\end{remark}
\begin{proof}
Let $f\in L^{\infty}_c$. For any
$g\in L^{p'(\cdot)}(\mathbb{R}^n)\subset L_{\rm loc}^1(\mathbb{R}^n)$
we have
\begin{equation}\label{le}
\int_{\mathbb{R}^n}|(Tf)(x)g(x)|dx
\le
c_n\int_{\mathbb{R}^n} Mf(x)Mg(x)dx.
\end{equation}
This inequality was proved in \cite[Theorem~3]{Le2}.
It follows easily by putting $Tf$ in place of $\varphi$ in
Theorem~\ref{th:Lerner} and by using Proposition~\ref{pr:relation} along
with Theorem~\ref{th:AlvarezPerez}. Notice that the application of
Theorem~\ref{th:Lerner} is justified due to the weak type $(1,1)$
of the operator $T$.

  From (\ref{le}), Lemma~\ref{le:Hoelder}, and the condition
$p,p'\in\mathcal{M}(\mathbb{R}^n)$ it follows that
$$
\int_{\mathbb{R}^n}|(Tf)(x)g(x)|dx
\le
c_nr_p\|Mf\|_{L^{p(\cdot)}}\|Mg\|_{L^{p'(\cdot)}}
\le
c_p\|f\|_{L^{p(\cdot)}}\|g\|_{L^{p'(\cdot)}},
$$
where $c_p:=c_nr_p\|M\|_{\mathcal{B}(L^{p(\cdot)})}
\|M\|_{\mathcal{B}(L^{p'(\cdot)})}$.
Then, by (\ref{eq:Orlicz-Luxemburg}),
$$
\|Tf\|_{L^{p(\cdot)}}
\le
\|Tf\|_{L^{p(\cdot)}}^0
\le
c_p\|f\|_{L^{p(\cdot)}}.
$$
  From the latter inequality and Lemma~\ref{le:density}
we get the theorem.
\end{proof}
\subsection{Zygmund spaces $L\log L(Q)$ and $L_{\exp}(Q)$}
\label{sec:Zygmund}
Let $Q$ be a cube in $\mathbb{R}^n$. The space
$L\log L(Q)$ consists of all measurable functions $f$ on $Q$ for which
$$
\int_Q|f(x)|\log^+|f(x)|dx<\infty
$$
(here $\log^+t:=\max\{\log t,0\}$). The space
$L_{\exp}(Q)$ consists of all measurable functions $f$ on $Q$ for which
there exists a $\lambda=\lambda(f)$ such that
$$
\int_Q\exp(\lambda|f(x)|)dx<\infty.
$$
It is well known (see, e.g., \cite[Chap.~4, Section~6]{BS})
that $L\log L(Q)$ and $L_{\exp}(Q)$ can be equip\-ped with the norms
\begin{eqnarray*}
\|f\|_{L\log L(Q)} &:=& \int_0^{|Q|}f^*(t)\log(|Q|/t)dt=\int_0^{|Q|}f^{**}(t)dt,
\\
\|f\|_{L_{\exp}(Q)} &:=&
\sup_{t\in(0,|Q|)}\frac{f^{**}(t)}{1+\log(|Q|/t)},
\end{eqnarray*}
respectively. Easy manipulations with $f^*$ and $f^{**}$ lead us
to the following well known H\"older inequality for Zygmund
spaces. If $f\in L\log L(Q)$ and $g\in L_{\exp}(Q)$, then $fg\in L^1(Q)$
and
\begin{equation}\label{holder}
\int_Q|f(x)g(x)|dx\le 2\|f\|_{L\log L(Q)}\|g\|_{L_{\exp}(Q)}.
\end{equation}
\subsection{A commuting relation for singular integrals}
The following statement represents one of the numerous variations
on the theme of adjoint operators (cf. \cite[Chap.~2, Section~5.3]{St1}),
and it seems to be known. We shall give its proof for the sake of
completeness.
\begin{prop}\label{commut}
Suppose $f$ and $\varphi$ are supported in a cube $Q\subset\mathbb{R}^n$.
If $\varphi\in L^{\infty}(Q)$ and $f\in L\log L(Q)$, then
\begin{equation}\label{comm}
\int_{\mathbb{R}^n}Tf(x)\varphi(x)dx
=
-\int_{\mathbb{R}^n}f(x)T\varphi(x)dx.
\end{equation}
\end{prop}
\begin{proof}
Set $K_{\varepsilon}(y):=K(y)\chi_{|y|>\varepsilon}$,
$T_{\varepsilon}f:=f*K_{\varepsilon}$, and
$$
T^*f(x):=\sup_{\varepsilon>0}|T_{\varepsilon}f(x)|.
$$
Since $K_{\varepsilon}$ is odd (because $K$ does), and the double
integral
$$
\int_{\mathbb{R}^n}\int_{\mathbb{R}^n}K_{\varepsilon}(x-y)f(y)\varphi(x)dydx
$$
converges absolutely, we clearly have
\begin{equation}\label{comm1}
\int_{\mathbb{R}^n} T_{\varepsilon}f(x)\varphi(x)dx
=
-\int_{\mathbb{R}^n}f(x)T_{\varepsilon}\varphi(x)dx\quad (\varepsilon>0).
\end{equation}
Next, $|(T_{\varepsilon}f)\varphi|\le |\varphi|(T^*f)$ and
$|f(T_{\varepsilon}\varphi)|\le |f|(T^*\varphi)$.

The conditions on $f$ and $\varphi$ imply $T^*f\in L^1(Q)$ and
$T^*\varphi\in L_{\exp}(Q)$, respectively (see, e.g.,
\cite[Chap.~2, Section~6.2]{St1}). Since $\varphi$ is supported in $Q$ and
$\varphi\in L^\infty(Q)$, we have $\varphi (T^*f)\in L^1(\mathbb{R}^n)$.
On the other hand, by the generalized H\"older inequality (\ref{holder}),
$f(T^*\varphi)\in L^1(\mathbb{R}^n)$.

Hence, letting $\varepsilon\to 0$ in (\ref{comm1}) and using the
dominated convergence theorem, we get (\ref{comm}).
\end{proof}
\subsection{Interpolation in Banach lattices}
We fix here some terminology (cf. \cite[Chap.~1]{CN03} and
\cite{Calderon64}). Let $(\mathcal{R},\Sigma,\mu)$ be a measure space and
$X$ be a Banach space of (equivalence classes of a.e. equal) real
valued measurable functions on $\mathcal{R}$ such that if $|g|\le|f|$
a.e., where $f\in X$ and $g$ is measurable, then $g\in X$ and
$\|g\|_X\le\|f\|_X$. The space $X$ is called a \textit{Banach
lattice} on $(\mathcal{R},\Sigma,\mu)$. The \textit{K\"othe dual} or
\textit{associate space} $X'$ of any Banach lattice $X$ on
$(\mathcal{R},\Sigma,\mu)$ is defined to be the space of real valued
measurable functions $g$ on $\mathcal{R}$ for which $fg\in
L^1(\mathcal{R},\Sigma,\mu)$ for each $f\in X$. For every $g\in X'$, put
\[
\|g\|_{X'}:=\sup\left\{\int_\mathcal{R}|fg|d\mu\ : \quad f\in X,\quad
\|f\|_X\le 1\right\}.
\]
To ensure that this is a norm rather than a seminorm we must
assume that $X$ is \textit{saturated}, that is, every $E\in\Sigma$
with $\mu(E)>0$ has a measurable subset $F$ of finite positive
measure for which $\chi_F\in X$.

Let $X_0$ and $X_1$ be Banach lattices on $(\mathcal{R},\Sigma,\mu)$ and
$0<\theta<1$. The \textit{Calder\'on product} (see
\cite[p.~123]{Calderon64}) consists of all real valued measurable
functions $f$ such that a.e. pointwise inequality
$|f|\le\lambda|f_0|^{1-\theta}|f_1|^\theta$ holds for some
$\lambda>0$ and elements $f_j$ in $X_j$ with $\|f_j\|_{X_j}\le 1$
for $j=0,1$. The norm of $f$ in $X_0^{1-\theta}X_1^\theta$ is
defined to be the infimum of all values $\lambda$ appearing in the
above inequality.    From results of \cite[Sections~6, 7, and~13.6]{Calderon64}
one can extract the following interpolation theorem.
\begin{theorem}\label{th:Calderon}
Let $X_0$ and $X_1$ be real Banach
lattices, one of which is reflexive. Let $A$ be a linear operator
bounded on $X_0$ and $X_1$. Then $A$ is bounded on
$X_0^{1-\theta}X_1^\theta$ and
\[
\|A\|_{\mathcal{B}(X_0^{1-\theta}X_1^\theta)}
\le 2
\|A\|_{\mathcal{B}(X_0)}^{1-\theta}
\|A\|_{\mathcal{B}(X_1)}^\theta.
\]
\end{theorem}
The following remarkable formula was proved by Lozanovski{\u\i}
\cite[Theorem~2]{Lozanovskii69} under some additional assumptions.
Cwikel and Nilsson relaxed assumptions on $X_0$ and $X_1$ and proved
the following (see \cite[Theorem~7.2]{CN03}).
\begin{theorem}\label{th:Lozanovskii}
For arbitrary Banach lattices $X_0$ and $X_1$,
\[
(X_0^{1-\theta}X_1^\theta)'=(X_0')^{1-\theta}(X_1')^\theta, \quad
0<\theta<1,
\]
with equality of the norms.
\end{theorem}

We refer to \cite[Chap.~15]{Maligranda89}
for generalizations of Theorems~\ref{th:Calderon}
and~\ref{th:Lozanovskii} to the case of so-called
\textit{Calder\'on-Lozanovski{\u\i} spaces}.

A Banach lattice $X$ is said to have the Fatou property if for
every a.e. pointwise increasing sequence $f_n$ of non-negative
functions in $X$ with $\sup\limits_n\|f_n\|_X<\infty$,
the function $f$, defined by
$f(x):=\lim\limits_{n\to\infty}f_n(x)$, is in $X$ and
$\|f\|_X=\lim\limits_{n\to\infty}\|f_n\|_X$. It is well known
(see, e.g., \cite[p.~452]{Zaanen67}) that if $X$ is a saturated
Banach lattice, then $X=X''$ isometrically if and only if $X$ has
the Fatou property.
\begin{cor}\label{co:Lozanovskii}
If $X$ is a saturated Banach lattice with the Fatou property and
$X'$ is its associate space, then
\begin{equation}\label{eq:Loz-1}
X^{1/2}(X')^{1/2}=L^2
\end{equation}
with equality of the norms.
\end{cor}
\begin{proof}
This result is contained in \cite[Theorem~5]{Lozanovskii69} in a
slightly different form. For the convenience of the readers we
reproduce here its proof from \cite[p.~185]{Maligranda89}.

Since $X$ is saturated, so is $X'$. Clearly, $X'$ has the Fatou
property. By the hypothesis, $X$ has the Fatou property too. Then
$X=X''$ and $X'=X'''$ with equalities for the norms. Put
$Z=X^{1/2}(X')^{1/2}$. Applying Theorem~\ref{th:Lozanovskii} with
$\theta=1/2$, we get
\begin{eqnarray}
Z' &=& (X')^{1/2}(X'')^{1/2}=(X')^{1/2}X^{1/2}=Z, \label{eq:Loz-2}
\\
Z'' &=& (X'')^{1/2}(X''')^{1/2}=X^{1/2}(X')^{1/2}=Z \nonumber
\end{eqnarray}
with equalities of the norms.

If $f\in Z=Z'=Z''$ is a non-zero function, then
\[
\|f\|_{Z'}=\|f\|_{Z''} 
= 
\sup_{\|g\|_{Z'}\le 1}\int_\mathcal{R}|fg|d\mu 
\ge
\int_\mathcal{R}\frac{f^2}{\|f\|_{Z'}}d\mu =
\frac{\|f\|_{L^2}^2}{\|f\|_{Z'}}.
\]
Hence,
\begin{equation}\label{eq:Loz-3}
\|f\|_Z=\|f\|_{Z'}=\|f\|_{Z''}\ge\|f\|_{L^2}.
\end{equation}
By duality, from the latter inequality we get
\begin{equation}\label{eq:Loz-4}
\|g\|_{L^2}\ge \|g\|_{Z'} \quad \mbox{for all}\quad g\in Z'.
\end{equation}
    From (\ref{eq:Loz-3}) and (\ref{eq:Loz-4}) it follows that
$L^2=Z'$ isometrically. Combining the latter equality with
(\ref{eq:Loz-2}), we arrive at (\ref{eq:Loz-1}).
\end{proof}

We shall apply the results of this section in the following form.
\begin{theorem}\label{th:interpolation}
Let $p:\mathbb{R}^n\to[1,\infty)$ be a measurable function
satisfying {\rm (\ref{eq:reflexivity})}. Let $A$ be a linear
operator bounded on $L^{p(\cdot)}(\mathbb{R}^n)$ and
$L^{p'(\cdot)}(\mathbb{R}^n)$. Then $A$ is bounded on $L^2(\mathbb{R}^n)$ and
\begin{equation}\label{eq:interpolation-1}
\|A\|_{\mathcal{B}(L^2)}
\le 2\sqrt{r_p}
\|A\|_{\mathcal{B}(L^{p(\cdot)})}^{1/2}
\|A\|_{\mathcal{B}(L^{p'(\cdot)})}^{1/2}.
\end{equation}
\end{theorem}
\begin{proof}
It is easy to see that $L^{p(\cdot)}(\mathbb{R}^n)$ is a saturated Banach
lattice on $\mathbb{R}^n$ equipped with the Lebesgue measure. Moreover,
$L^{p(\cdot)}(\mathbb{R}^n)$ has the Fatou property (see, e.g.,
\cite[Proposition~1.3]{ELN99}). So, we can apply Theorem~\ref{th:Calderon}
and Corollary~\ref{co:Lozanovskii}.

   From the obvious equality $r_p=r_{p'}$ and (\ref{eq:Orlicz-Luxemburg}) we get
\begin{equation}\label{eq:interpolation-2}
\|A\|_{\mathcal{B}([L^{p(\cdot)}]')}
\le r_p
\|A\|_{\mathcal{B}(L^{p'(\cdot)})}.
\end{equation}
Applying Theorem~\ref{th:Calderon} with $\theta=1/2$ and
Corollary~\ref{co:Lozanovskii} to $X_0=X=L^{p(\cdot)}(\mathbb{R}^n)$ and
$X_1=X'=[L^{p(\cdot)}(\mathbb{R}^n)]'$, we get
\begin{equation}\label{eq:interpolation-3}
\|A\|_{\mathcal{B}(L^2)}
=
\|A\|_{\mathcal{B}\big((L^{p(\cdot)})^{1/2}([L^{p(\cdot)}]')^{1/2}\big)}
\le  2
\|A\|_{\mathcal{B}(L^{p(\cdot)})}^{1/2}
\|A\|_{\mathcal{B}([L^{p(\cdot)}]')}^{1/2}.
\end{equation}
Combining (\ref{eq:interpolation-2}) and
(\ref{eq:interpolation-3}), we arrive at
(\ref{eq:interpolation-1}).
\end{proof}
\section{Proof of Theorem \ref{Main}}
\label{sec:proof} We start with Part (a). The proof of this part
is similar to the proof of Theorem~\ref{th:CZ}.

Let $f\in L_c^\infty$ and
$g\in L^{p'(\cdot)}(\mathbb{R}^n)\subset L_{\rm loc}^1(\mathbb{R}^n)$.
In view of \cite[Theorem~1.1]{Perez}, the function $[b,T]f$
satisfies (\ref{eq:Lerner}). Thus, putting $[b,T]f$ in place of
$\varphi$ in Theorem~\ref{th:Lerner} and then applying
Proposition~\ref{pr:relation} and Theorem~\ref{th:Perez}, we get
\begin{eqnarray*}
\int_{\mathbb{R}^n}\Big|\big([b,T]f\big)(x)g(x)\Big|dx
&\le&
c_n(1/\lambda_n)^{1/\delta}
\int_{\mathbb{R}^n}([b,T]f)_\delta^\#(x)Mg(x)dx
\\
&\le&
c'\|b\|_*\int_{\mathbb{R}^n}M(Tf)(x)Mg(x)dx
\\
&+&
c'\|b\|_*\int_{\mathbb{R}^n}MMf(x)Mg(x)dx
\end{eqnarray*}
with $c':=c_nc_{\delta,n}(1/\lambda_n)^{1/\delta}$.
  From the latter inequality, Lemma~\ref{le:Hoelder}, Theorem~\ref{th:CZ},
and the condition $p,p'\in\mathcal{M}(\mathbb{R}^n)$ it follows that
\begin{eqnarray}
\label{eq:Main-1}
&&
\int_{\mathbb{R}^n}\Big|\big([b,T]f\big)(x)g(x)\Big|dx
\\
\nonumber
&&\le
c'r_p\|b\|_*
\Big(\|M(Tf)\|_{L^{p(\cdot)}}+\|MMf\|_{L^{p(\cdot)}}\Big)
\|Mg\|_{L^{p'(\cdot)}}
\\
\nonumber
&&\le
C_p\|b\|_*\|f\|_{L^{p(\cdot)}}\|g\|_{L^{p'(\cdot)}},
\end{eqnarray}
where
$$
C_p:=
c'r_p
\|M\|_{\mathcal{B}(L^{p(\cdot)})}
\|M\|_{\mathcal{B}(L^{p'(\cdot)})}
\Big(
\|T\|_{\mathcal{B}(L^{p(\cdot)})}
+
\|M\|_{\mathcal{B}(L^{p(\cdot)})}
\Big).
$$
  From (\ref{eq:Main-1}) and (\ref{eq:Orlicz-Luxemburg})
we obtain for all $f\in L_c^\infty$,
$$
\|[b,T]f\|_{L^{p(\cdot)}}
\le
\|[b,T]f\|_{L^{p(\cdot)}}^0
\le
C_p\|b\|_{*}\|f\|_{L^{p(\cdot)}}.
$$
By Lemma \ref{le:density}, $[b,T]$ can be extended by continuity
to a bounded linear operator on $L^{p(\cdot)}(\mathbb{R}^n)$, and
$\|[b,T]\|_{\mathcal{B}(L^{p(\cdot)})}\le C_p\|b\|_*$.
Part (a) is proved.

We turn now to the proof of Part (b). Suppose $b\in L\log L(Q)$
for any cube $Q$ and $[b,T]$ is bounded on
$L^{p(\cdot)}(\mathbb{R}^n)$. Let $f\in L_c^\infty$ and
$\varphi\in L^{p(\cdot)}(\mathbb{R}^n)$. For natural $k$ set
$\varphi_k:=\min\{|\varphi|,k\}\chi_{B(0,k)}$, where $B(0,k)$
is the ball of radius $k$ centered at the origin. Let also
$\varphi'_k=\varphi_k\,\text{sgn}\,([b,T]f)$.

Clearly, $bf$ and $b\varphi'_k$ belong to $L\log L(Q)$ for any
cube $Q\subset\mathbb{R}^n$ and any $k$. Applying Proposition
\ref{commut} yields
\begin{eqnarray*}
\int_{\mathbb{R}^n}|\big([b,T]f\big)(x)\varphi_k(x)|dx&=&
\int_{\mathbb{R}^n}\big([b,T]f\big)(x)\varphi'_k(x)dx \\
&=&\int_{\mathbb{R}^n}\big([b,T]\varphi'_k\big)(x)f(x)dx.
\end{eqnarray*}
Hence, by Lemma~\ref{le:Hoelder},
\begin{eqnarray*}
\int_{\mathbb{R}^n}|\big([b,T]f\big)(x)\varphi_k(x)|dx &\le&
r_p\|[b,T]\|_{\mathcal{B}(L^{p(\cdot)})}
\|\varphi'_k\|_{L^{p(\cdot)}} \|f\|_{L^{p'(\cdot)}}
\\
&\le&
r_p\|[b,T]\|_{\mathcal{B}(L^{p(\cdot)})}
\|\varphi\|_{L^{p(\cdot)}}
\|f\|_{L^{p'(\cdot)}}.
\end{eqnarray*}
By the Fatou convergence theorem,
$$
\int_{\mathbb{R}^n}|\big([b,T]f\big)(x)\varphi(x)|dx \le
r_p\|[b,T]\|_{\mathcal{B}(L^{p(\cdot)})}
\|\varphi\|_{L^{p(\cdot)}} \|f\|_{L^{p'(\cdot)}}.
$$
Hence, taking into account (\ref{eq:Orlicz-Luxemburg}), for every $f\in L_c^\infty$,
$$
\|[b,T]f\|_{L^{p'(\cdot)}}
\le
\|[b,T]f\|_{L^{p'(\cdot)}}^0\le
r_p\|[b,T]\|_{\mathcal{B}(L^{p(\cdot)})}
\|f\|_{L^{p'(\cdot)}}.
$$
  From the latter inequality and Lemma~\ref{le:density} we deduce that $[b,T]$
is bounded on $L^{p'(\cdot)}(\mathbb{R}^n)$ and
$\|[b,T]\|_{\mathcal{B}(L^{p'(\cdot)})}
\le
r_p\|[b,T]\|_{\mathcal{B}(L^{p(\cdot)})}$.
In view of Theorem~\ref{th:interpolation},
it follows that $[b,T]$ is bounded on $L^2(\mathbb{R}^n)$ and
\begin{equation}\label{eq:nec-1}
\|[b,T]\|_{\mathcal{B}(L^2)}\le
2r_p
\|[b,T]\|_{\mathcal{B}(L^{p(\cdot)})}.
\end{equation}
On the other hand, by Janson's theorem \cite{Janson}, if $[b,T]$
is bounded on $L^2(\mathbb{R}^n)$, then $b\in BMO(\mathbb{R}^n)$.
Moreover, from the proof in \cite{Janson} one can see that
there exists a positive constant $c_2(K)$ such that
\begin{equation}\label{eq:nec-2}
\|b\|_*\le c_2(K)\|[b,T]\|_{\mathcal{B}(L^2)}.
\end{equation}
Combining (\ref{eq:nec-1}) and (\ref{eq:nec-2}), we arrive at
Part (b) with $C_p':=2r_pc_2(K)$.
Theorem~\ref{Main} is proved.
\section{Concluding remarks}
\label{sec:concluding}
\begin{remark}\label{rem:final}
We have learned recently that Diening has obtained a new characterization
of the class $\mathcal{M}(\mathbb{R}^n)$. In particular,
$p\in \mathcal{M}(\mathbb{R}^n)$ if and only if
$p'\in \mathcal{M}(\mathbb{R}^n)$ \cite[Theorem~8.1]{Di}
and $p\in\mathcal{M}(\mathbb{R}^n)$ implies $(p/s)'\in\mathcal{M}(\mathbb{R}^n)$
for some $s\in(0,1)$ \cite[Corollary~8.8]{Di}. So, the condition
$p'\in\mathcal{M}(\mathbb{R}^n)$ in Theorems~\ref{Main} and~\ref{th:CZ}
can be removed. For singular integral operators this is noticed
already in \cite[Theorem~8.14]{Di}.
\end{remark}

\begin{remark}
Part (a) of Theorem~\ref{Main} holds for a more general class of
Calder\'on-Zygmund operators \cite{Jo}.
\end{remark}

\begin{remark}
All results of this paper can be extended to the context of Banach function
spaces in the sense of Luxemburg (see \cite[Chap.~1, Definition~1.1]{BS}).
   From \cite[Chap.~1, Theorem~3.11 and Corollary~4.4]{BS} it follows that
if $X(\mathbb{R}^n)$ is a reflexive Banach function space, then
$L_c^\infty$ is dense in $X(\mathbb{R}^n)$ and in its associate space
$X'(\mathbb{R}^n)$. Theorems~\ref{Main} and~\ref{th:CZ} remain valid
for reflexive Banach function spaces under the assumption that the
Hardy-Littlewood maximal function is bounded on $X(\mathbb{R}^n)$
and on $X'(\mathbb{R}^n)$. The proofs of these statements are minor
modifications of those for Theorems~\ref{Main} and~\ref{th:CZ}.
\end{remark}
\section*{Acknowledgments}
The first author is supported by F.C.T. (Portugal) grant
SFRH/ BPD/11619/2002. The second author is grateful for the
support given by ``Centro de Matem\'a\-tica e Aplica\c{c}\~oes'',
Instituto Superior T\'ecnico, Lisbon, Portugal, during his short
stay  in Lisbon in January--February of 2004.

We would like to thank Lars Diening and Ale\v{s}
Nekvinda for sharing with us their preprints \cite{Diening02,Di,Nekvinda02}.

\end{document}